\documentclass[11pt,reqno]{amsart}

\theoremstyle{plain}
\newtheorem{theorem}{Theorem} [section]
\newtheorem{corollary}[theorem]{Corollary}
\newtheorem{lemma}[theorem]{Lemma}
\newtheorem{proposition}[theorem]{Proposition}

\theoremstyle{definition}
\newtheorem{definition}[theorem]{Definition}

\newtheorem{remark}[theorem]{Remark}

\def\N{{\mathbb{N}}}
\def\R{{\mathbb{R}}}
\def\Z{{\mathbb{Z}}}

\def\F{{\mathcal{F}}}
\def\clspan{{\overline{\mathrm{span}}}}
\def\comp{{\mathrm{C}}}
\def\fhat{{\widehat{f}}}
\def\ghat{{\widehat{g}}}

\def\Jc{{\mathcal{J}}}
\def\Sf{{\mathfrak S}}

\def\Tc{{\mathcal{T}}}
\def\phihat{{\widehat{\varphi}}}

\newcommand{\bigabs}[1]{{\bigl|#1\bigr|}}

\newcommand{\CHI}{\hbox{\raise .4ex \hbox{$\chi$}}}
\newcommand{\chiBk}{{\CHI_{B_k}}}
\newcommand{\medcup}{\operatornamewithlimits{\textstyle\bigcup}}
\newcommand{\Dim}{{\mathrm{dim}}}
\newcommand{\ip}[2]{\langle#1,#2\rangle}
\newcommand{\bigip}[2]{\bigl\langle #1, \, #2 \bigr\rangle}
\newcommand{\Bigip}[2]{\Bigl\langle #1, \, #2 \Bigr\rangle}
\newcommand{\norm}[1]{\lVert#1\rVert}

\newcommand{\bigparen}[1]{\bigl(#1\bigr)}
\newcommand{\Bigparen}[1]{\Bigl(#1\Bigr)}
\newcommand{\biggparen}[1]{\biggl(#1\biggr)}
\newcommand{\rank}{{\mathrm{rank}}}
\newcommand{\set}[1]{\{#1\}}
\newcommand{\bigset}[1]{\bigl\{#1\bigr\}}
\newcommand{\Bigset}[1]{\Bigl\{#1\Bigr\}}
\newcommand{\Span}{\mathrm{span}}
\newcommand{\supp}{{\mathrm{supp}}}

\newcommand{\dotoplus}{\operatornamewithlimits{\dot{\oplus}}}
\newcommand{\dotbigoplus}{\operatornamewithlimits{\dot{\bigoplus}}}

\def\subset{\subseteq}

\def\EQ{\, = \,}

\def\GE{\, \ge \,}

\def\LE{\, \le \,}
\def\IN{\, \in \,}
\def\TO{\, \to \,}

\def\IMPLIES{\quad \implies \quad}

\begin{document}

\title{Invariance of a Shift-Invariant Space}
\author[A. Aldroubi, C. Cabrelli, C. Heil, K. Kornelson,
and U. Molter]{Akram Aldroubi, Carlos Cabrelli, Christopher Heil, \\
Keri Kornelson, and Ursula Molter}

\address{\textrm{(A. Aldroubi)}
Department of Mathematics,
Vanderbilt University,
Nashville, Tennessee 37240-0001 USA}
\email{aldroubi@math.vanderbilt.edu}

\address{\textrm{(C. Cabrelli)}
Departamento de Matem\'atica,
Facultad de Ciencias Exac\-tas y Naturales,
Universidad de Buenos Aires, Ciudad Universitaria, Pabell\'on I,
1428 Buenos Aires, Argentina and
CONICET, Consejo Nacional de Investigaciones
Cient\'ificas y T\'ecnicas, Argentina}
\email{cabrelli@dm.uba.ar}

\address{\textrm{(C. Heil)}
School of Mathematics,
Georgia Institute of Technology,
Atlanta, Georgia 30332-0160 USA}
\email{heil@math.gatech.edu}

\address{\textrm{(K. Kornelson)}
Department of Mathematics,
 University of Oklahoma, 
 Norman, Oklahoma 73019-0315 USA
}
\email{kkornelson@math.ou.edu}

\address{\textrm{(U. Molter)}
Departamento de Matem\'atica,
Facultad de Ciencias Exac\-tas y Naturales,
Universidad de Buenos Aires, Ciudad Universitaria, Pabell\'on I,
1428 Buenos Aires, Argentina and
CONICET, Consejo Nacional de Investigaciones
Cient\'ificas y T\'ecnicas, Argentina}
\email{umolter@dm.uba.ar}

\thanks{The research of
A.~Aldroubi was supported in part by NSF Grant DMS-0807464.
The research of C.~Cabrelli and U.~Molter was partially supported by
Grants FONCyT PICT 2006-00177, CONICET, PIP 5650, UBACyT X058 and X108.
The research of C.~Heil was supported in part by NSF  Grant DMS-0806532.
The research of K.~Kornelson was supported in part by NSF Grant DMS-0701164, Grinnell College CSFS, the Woodrow Wilson Fellowship Foundation, and the IMMERSE program at the University of Nebraska-Lincoln.}

\subjclass[2000]{Primary 42C40; Secondary 42C15, 46C99}

\keywords{Dimension function, Fiber space, frames,
Gramian operator, Range function, 
Riesz basis,
shift-invariant space,
translation-invariant space}

\date{\today}

\begin{abstract}
A shift-invariant space is a space of functions that is
invariant under integer translations.
Such spaces are often used as models for spaces of signals and images
in mathematical and engineering applications.
This paper characterizes those shift-invariant subspaces~$S$
that are also invariant under additional (non-integer) translations.
For the case of finitely generated spaces, these spaces are characterized
in terms of the generators of the space.
As a consequence, it is shown that principal shift-invariant spaces
with a compactly supported generator cannot be invariant under any
non-integer translations.
\end{abstract}

\maketitle

\section{Introduction}

A \emph{shift-invariant space} (SIS) is a space of functions that is
invariant under integer translations.

They have applications throughout mathematics and engineering, as such
spaces are often used as models for spaces of signals and images, see
\cite{Gro01}, \cite{HW96}, \cite{Mal98}.

One example of a shift-invariant space is the Paley--Wiener space of
functions that are bandlimited to $[-1/2,1/2]$:
$$PW(\R)
\EQ \bigset{f \in L^2(\R) : \supp(\fhat\,) \subseteq [-\tfrac12,\tfrac12]}.$$
This SIS has the property that it is not only invariant under integer
translations, but it is in fact invariant under every real translation.
A space with this property is said to be \emph{translation-invariant}.
A classical theorem of Fourier analysis (often attributed to Wiener, see for example \cite{Hel64}),  completely characterizes the closed
translation-invariant subspaces of $L^2(\R)$ as being of the form
$$\set{f \in L^2(\R) : \supp(\fhat\,) \subseteq A}$$
where $A \subseteq \R$ is measurable.

In many applications, it is desirable to have a shift-invariant
space that possesses extra invariances \cite{CS03}, \cite{Web00}.
In this paper we characterize those shift-invar\-iant subspaces~$S$
that are not only invariant under integer translations,
but are also invariant under some particular set of
translations~$M \subseteq \R$.
We show that there are only two possibilities:
\begin{itemize}
\item
either $S$ is translation-invariant, or

\item
there exists an $n \in \N$ such that~$S$ is invariant under
translations by multiples of $1/n$, but not invariant under
translations by $1/m$ with $m>n$.
\end{itemize}

\noindent
We give several characterizations of those shift-invariant
spaces that are $\frac1n\Z$-invariant.
A trivial way to create such a space is to fix a
function $g \in L^2(\R)$ and set
$$S \EQ \clspan\bigset{g(x - \tfrac{k}{n}) : k \in \Z},$$
the closed span of the $\frac1n\Z$ translates of~$g$.
However, we are interested in the more
subtle question of recognizing when a given SIS is $\frac1n\Z$-invariant.
For example, in many applications one is presented with a SIS of the form
$$S \EQ \clspan\bigset{g(x-k) : k \in \Z},$$
and it is not obvious whether such a space possesses any invariants
other than translation by integers.
We completely determine the invariances of such a space in terms of
properties of~$g$ and, more generally, characterize any SIS
that is $\frac1n\Z$-invariant.

One interesting corollary of our characterization is that the
shift-invariant space generated by a compactly supported function
is not invariant under any translations other than~$\Z$.
Thus, the shift-invariant spaces associated with compactly supported
multiresolution analyses and wavelets are already ``maximally invariant.''

\section{Notation and Definitions}

We normalize the Fourier transform of $f \in L^1(\R)$ as
$$\fhat(\omega) \EQ \int_{-\infty}^\infty f(x) \, e^{-2\pi i\omega x} dx.$$
The Fourier transform extends to a unitary operator on $L^2(\R)$.
Given $\F \subseteq L^2(\R)$, we set
$\widehat \F = \set{\fhat : f \in \F}$.

The translation operator $T_a$ is $T_a f(x) = f(x-a)$.
Note that $(T_a f)^\wedge(\omega) = e^{-2\pi i a \omega} \fhat(\omega)$.

A function $f$ is $b\Z$-\emph{periodic} if
$T_{bk}f = f$ for all $k \in \Z$.
A set $A \subseteq \R$ is $b\Z$-\emph{periodic} if its
characteristic function is $b\Z$-periodic.

A \emph{shift-invariant space} (SIS) is a closed subspace $S$ of $L^2(\R)$
that is invariant under integer translations.
We say that $S$ is $b\Z$-invariant if it is invariant under translation
by $bk$ for all $k \in \Z$.

Given $\F \subseteq L^2(\R)$, we define
$$\Tc_\Z(\F)
\EQ \set{T_j f: f \in \F, \, j \in \Z}.$$
The SIS generated by $\F$ is
$$\Sf(\F)
\EQ \clspan\bigparen{\Tc_\Z(\F)}
\EQ \clspan\set{T_j f: f \in \F, \, j \in \Z}.$$
We call $\F$ a \emph{set of generators} for $\Sf(\F)$.
When $\F = \set{f}$ consists of a single function, we simply write $\Sf(f)$.

The \emph{length} of a SIS $S$ is the minimum cardinality of the
sets $\F$ such that $S = \Sf(\F)$.
A SIS of length one is called a \emph{principal} SIS.
A SIS of finite length is a \emph{finitely generated} SIS.

We will write $W = U \dotoplus V$ to denote the \emph{orthogonal} direct
sum of closed subspaces of $L^2(\R)$,
i.e., the subspaces $U$, $V$ must be closed and orthogonal, and~$W$
is their direct sum.

The Lebesgue measure of a set $E \subseteq \R$ is denoted by $|E|$.

The cardinality of a finite set $F$ is denoted by $\#F$.

\section{Order of Invariance}

Let $S$ be a SIS.
If $\theta$ is a real number,  we will say that $S$ is {\it invariant under translations
by~$\theta$} or that $S$ is {\it $\theta$-invariant} if
$$f \in S \IMPLIES T_\theta f \in S.$$

We have the following Proposition.
\begin{proposition}\label{group}
Let $S$ be a SIS and define,
$$M = \{\theta \in \R: S \text{ is } \theta\text{-invariant} \}.$$
Then $M$ is a closed additive subgroup of $\R$ containing $\Z.$ 
\end{proposition}

\begin{proof}
Note that $\Z \subseteq M$ since $S$ is shift-invariant.
To see that M is closed, let $\{\theta_j\}$ be a sequence in $M$ such that
$\theta_j \rightarrow \theta.$
Then, given any $f \in S$, we have,
$$\norm{T_{\theta_j} f - T_{\theta} f}_2^2
\EQ \int_{-\infty}^\infty |f(x-\theta_j) - f(x-\theta)|^2 \, dx
\TO 0 \quad\text{as } j \to \infty.$$
So, since $S$ is closed and $T_{\theta_j} \in S$,  $T_{\theta} f$ must be in $S$ and therefore $\theta \in M.$

Let us prove now that M is indeed an additive subgroup of $\R$.
Clearly, M is closed under addition. Furthermore, if $ n,m \in \Z $ with $n > 0$ and $\theta \in M $ then  $ n\theta + m \in M.$

We need to see that $-\theta$ is in $M$ for each $\theta \in M.$
For this, let us first consider the case that $\theta$ is rational. So we can assume that $\theta = p/q$  with $q>0.$ 

Then, if $p/q \in M$ we have
$$ -\frac{p}{q} \EQ (q-1)\frac{p}{q} -p \in M.$$

Now, if $\theta \in M$ is irrational, then $D\equiv \{n\theta+m: n,m\in \Z, n > 0 \} \subset M.$
Since $D$ is dense in $\R$ and $M$ is closed, then $M = \R$ and so $-\theta$ is in $M$.
\end{proof}

Since  the only closed additive subgroups of $\R$ containing $\Z$ are $\frac{1}{n}\Z$ for some positive integer $n$  or the entire group $\R,$  we have the following.

\begin{proposition} \label{first}
Let $S$ be a SIS.
Then either $S$ is translation-invariant, or there exists
a maximum positive integer $n$ such that $S$ is $\frac1n\Z$-invariant.
\end{proposition}
 
Proposition~\ref{first} suggests the following definition.

\begin{definition}
Given a shift-invariant space $S$,
we say that $S$ has {\it invariance order~$n$} if~$n$ is the maximum positive
integer such that $S$ is $\frac1n\Z$-invariant.
If this maximum does not exist, we say that $S$ has {\it invariance
order~$\infty$}; in this case~$S$ is translation-invariant.
\end{definition}

\begin{remark}
Note that the invariance order of any SIS is at least~$1$,
since $S$ is $\Z$-invariant.
Also, if $S$ has invariance order $n$, then~$S$ is \emph{not}
invariant under translation by any real number $y$ in the range
$0 < y < 1/n$.
Furthermore, if $y \geq 1/n$, $S$ can \emph{only} be invariant under
translation by~$y$ if~$y$ is a multiple of~$1/d$ where~$d$ divides~$n$.
In particular, if the order of invariance of $S$ is a prime number $p$,
there exist no other integers $m > 1$ such that~$S$ is invariant
under translations by $1/m$.
\end{remark}

\section{Characterization of $\frac1n$-Invariance}

In this part we will characterize those shift-invariant spaces
that are $\frac1n\Z$-invariant.

\smallskip
\subsection{Notation}

We will use the following notation throughout the remainder of this paper.

Given a fixed positive integer~$n$, we partition the real line into~$n$ sets,
each of which is $n\Z$-periodic, as follows.
For $k=0,\dots,n-1$ define,
$$B_k \EQ \medcup_{j \in \Z} \,\, ([k,k+1)+nj).$$
Note that $B_k$ implicitly depends on the choice of~$n$.

Given a SIS $S \subseteq L^2(\R)$, we associate the following subspaces:
\begin{equation} \label{Sk}
U_k \EQ \set{f \in L^2(\R) :
             \fhat = \ghat \, \chiBk \text{ for some } g \in S },
\qquad k = 0,\dots,n-1.
\end{equation}
The spaces $U_k$ are mutually orthogonal since the sets~$B_k$
are disjoint (up to sets of measure zero).

If $f \in S$ and $0 \le k \le n-1$, then we let $f^k$ denote the
function defined by
$$\widehat{f^k} \EQ \fhat \, \chiBk.$$
Letting $P_k$ denote the orthogonal projection onto $\{ f : {\text supp}(\hat f) \subset B_k\}$, we have that
$$U_k = P_k(S) \text{ and } f^k \EQ P_k f.$$
Note that integer translations commute with the projections~$P_k$:
if $j \in \Z$ and $k = 0, \dots, n-1$, then
$$T_j P_k \EQ P_k T_j.$$

\smallskip
\subsection{Preliminary results}

We will need the following  result from \cite{dBDR94a}.

\begin{proposition}[\cite{dBDR94a}] \label{BDR2}
Let $f \in L^2(\R)$ be given.
If $g \in \Sf(f)$, then there exists a $\Z$-periodic function $m$
such that $\ghat = m \fhat$.

Conversely, if $m$ is a $\Z$-periodic function such that
$m \fhat \in L^2(\R)$,
then the function $g$ defined by $\ghat = m \fhat$ belongs to $\Sf(f)$.
\end{proposition}

We will also need a version of the preceding result for spaces that
are $\frac1n\Z$-invariant instead of shift-invariant.
This follows easily by rescaling.

\begin{corollary} \label{BDR2_cor}
Let $f \in L^2(\R)$ and $n \in \N$ be given, and set
$$\Sf(f,\tfrac1n\Z)
\EQ \clspan\set{T_{j/n} f:  \, j \in \Z}.$$
If $g \in \Sf(f,\tfrac1n\Z)$, then there exists a $n\Z$-periodic function~$m$
such that $\ghat = m \fhat$.

Conversely, if $m$ is an $n\Z$-periodic function such that
$m \fhat \in L^2(\R)$,
then the function $g$ defined by $\ghat = m \fhat$
belongs to $\Sf(f,\tfrac1n\Z)$.
\end{corollary}

\smallskip
\subsection{Characterization of $\frac1n$-invariance in terms of subspaces}

The periodicity of the $B_k$ sets yields the following lemma.

\begin{lemma} \label{lemma}
Let  $S$ be  a SIS. Assume that the subspace $U_k \subset S$. Then for each $k=0, \dots, n-1$,
 $U_k$ is a SIS that is also $\frac1n\Z$-invariant.
\end{lemma}

\begin{proof}
Fix $0 \le k \le n-1$, and choose any $f \in U_k$.
There exists a $g \in S$ such that $\fhat = \ghat \, \chiBk$.
Since $S$ is shift-invariant and $g \in S$, we have that
$e^{-2\pi i s\omega} \ghat(\omega)$ is in $\widehat S,$\ for all $s\in \Z$.
Hence
$$e^{-2\pi i s\omega} \, \fhat(\omega)
\EQ e^{-2\pi i s\omega} \, \ghat(\omega) \, \chiBk(\omega)
\IN \widehat U_k.$$
Therefore $T_s f \in U_k$, so $U_k$ is invariant under integer translates.

Suppose now that $f_j \in U_k$ and $f_j \to f$ in $L^2(\R)$.
Since $U_k \subseteq S$ and $S$ is closed, $f$ must be in $S$.
Further,
$$\norm{\fhat_j - \fhat \,}_2^2
\EQ \norm{(\fhat_j - \fhat\,) \, \chiBk}_2^2 +
    \norm{(\fhat_j - \fhat\,) \, \CHI_{B_k^\comp}}_2^2
\EQ \norm{\fhat_j - \fhat \, \chiBk}_2^2 +
    \norm{\fhat \, \CHI_{B_k^\comp}}_2^2.$$
Since the left-hand side converges to zero,
we must have that $\fhat \, \CHI_{B_k^\comp} = 0$~a.e.,
and that $\fhat_j \to \fhat \, \chiBk$ in $L^2(\R)$.
Since we also have $\fhat_j \to \fhat$, we conclude that
$$\fhat \EQ \fhat \, \chiBk \textrm{ a.e.}$$
Consequently $f \in U_k$, so $U_k$ is closed.

Finally, to see that $U_k$ is $\frac1n\Z$-invariant, define
$$h(\omega)
\EQ e^{-\frac{2\pi i \omega}{n}}
    \sum_{j=-k}^{n-1-k} e^{\frac{2 \pi i j}{n}} \, \CHI_{B_{k+j}}(\omega).$$ 
Note that $|h(\omega)|=1$ and that $h$ is $\Z$-periodic.
Furthermore, if $\omega \in B_k$ and $-k \le j \le n-1-k$, then
$\CHI_{B_{k+j}}(\omega)$ can be nonzero only when $j=0$.
Hence:
$$\omega \in B_k
\IMPLIES h(\omega) = e^{-\frac{2\pi i \omega}{n}}.$$
If $f \in U_k$ then, since $\supp(\hat f) \subseteq B_k$, we have
$$e^{-\frac{2 \pi i \omega}{n}} \fhat(\omega)
\EQ h(\omega) \, \fhat(\omega).$$
However, since $U_k$ is $\Z$-invariant, we have
by Proposition~\ref{BDR2} that $h \fhat \in \widehat{U_k}$.
Therefore $e^{-\frac{2 \pi i \omega}{n}} \fhat(\omega) \in \widehat{U_k}$,
which implies that $T_{1/n}f \in U_k$.
\end{proof}

This leads to the following characterization.

\begin{theorem} \label{main}
If $S \subseteq L^2(\R)$ is a SIS, then the following are equivalent.

\begin{enumerate}
\item[(a)] $S$ is $\frac1n\Z$-invariant.

\smallskip
\item[(b)] $U_k \subseteq S$ for $k=0, \dots, n-1$. \label{part-2}

\smallskip
\item[(c)] If $f \in S$, then $f^k = P_k f \in S$
for each $k = 0, \dots, n-1$.

\end{enumerate}

\noindent
Moreover, in case these hold we have that $S$ is the orthogonal direct sum
$$S \EQ U_0 \dotoplus \dots \dotoplus U_{n-1},$$
with each $U_k$ being a (possibly trivial) $\frac1n\Z$-invariant SIS.
\end{theorem}
\begin{proof}
(a) $\Rightarrow$ (b).
Assume that $S$ is $\frac1n\Z$-invariant and fix $0 \le k \le n-1$
and $f \in U_k$.
By definition of $U_k$, we have that
$\fhat = \ghat \, \chiBk$ for some $g \in S$.
Since $\chiBk$ is $n\Z$-periodic and bounded,
Corollary~\ref{BDR2_cor} implies that
$f \in \Sf(g,\tfrac1n\Z) \subseteq S$.

\medskip
(b) $\Rightarrow$ (a).
Suppose that $U_k \subseteq S$ for each $k = 0,\dots,n-1$.

Note that Lemma~\ref{lemma} implies that $U_k$ is $\frac1n\Z$-invariant,
and we also have that the $U_k$ are mutually orthogonal
since the sets $B_k$ are disjoint.

Suppose that $f \in S$.
Then $f = f^0 + \dots + f^{n-1}$
where
$\widehat{f^k} = \fhat \, \chiBk$.
This implies that
$f \in U_0 \dotoplus \dots \dotoplus U_{n-1}$,
and consequently~$S$ is the orthogonal direct sum
$$S \EQ U_0 \dotoplus \dots \dotoplus U_{n-1}.$$

As each $U_k$ is $\frac1n\Z$-invariant, it follows that~$S$
is $\frac1n\Z$-invariant as well.

\medskip
(b) $\Leftrightarrow$ (c).
This is a restatement of the definition of $U_k$.
\end{proof}

\begin{corollary} \label{colo-II}
Let $S$ be a SIS.
If there exists a $k \in \set{0,\dots,n-1}$ such that 
$\supp(\fhat\,) \subseteq B_k$ for all $f \in S$,
then $S$ is $\frac1n\Z$-invariant.
\end{corollary}

\begin{remark}
It is interesting  to note that the subspaces $U_k$ satisfy:

$$
U_k = P_k (S) = \{ f \in L^2(\R): {\text supp}(\hat f) \subset B_k\} \cap S .
$$  

That is, the projections and the restrictions of $S$ yield  valid tests for $\frac{1}{n}\Z$-invariance.
\end{remark}

\smallskip
\subsection{Characterization of $\frac1n$-invariance in terms of generators}

We will show now that the conditions for $\frac1n\Z$-invariance
can be formulated in terms of properties of a set of generators of the SIS.

\begin{theorem}
Let $\F$ be a set of generators for a SIS $S$, i.e., $S=\Sf(\F)$.
Then the following statements are equivalent.

\begin{enumerate}
\item[(a)] $S$ is $\frac1n\Z$-invariant.

\smallskip
\item[(b)]
$P_k\F = \set{f^k: f \in \F} \subseteq S$ for $k=0, \dots, n-1$.
\end{enumerate}
\end{theorem}

\begin{proof}
(a) $\Rightarrow$ (b).
This is a consequence of Theorem~\ref{main}.

\medskip
(b) $\Rightarrow$ (a).
Suppose that statement~(b) holds.
Then, by hypothesis, $V_k = \Sf(P_k\F) \subseteq S$,
and, by Corollary~\ref{colo-II}, $V_k$ is $\frac1n\Z$-invariant.
Furthermore, $V_j \perp V_k$ when $j \ne k$.
If $f \in \F$, then
$f = f^0 + \cdots + f^{n-1} \in V_1 \dotoplus \cdots \dotoplus V_{n-1}$.
Consequently,
$S = \Sf(\F) = V_1 \dotoplus \cdots \dotoplus V_{n-1}$.
As each $V_k$ is $\frac1n\Z$-invariant, it follows that $S$ is as well.
\end{proof}

It is known that it is always possible to choose a (possibly infinite)
set of generators of a SIS in such a way that the integer translates of
the generators actually forms a frame for the SIS (see Theorem \ref{Bownik}).
This is particularly important in applications, and
we examine this situation next.

Recall that a countable collection of vectors 
$\set{v_{\alpha} : \alpha \in \Lambda}$
forms a \emph{frame} for a Hilbert space~$H$ if there 
exist constants $A$, $B$ (called \emph{frame bounds}) such that
\begin{equation} \label{frame-condition}
\forall\, w \in H, \quad
A \, \norm{w}^2
\LE \sum_{\alpha \in \Lambda} |\ip{w}{v_\alpha}|^2
\LE B \norm{w}^2.
\end{equation}
If we can take $A = B = 1$, then the frame is called a \emph{Parseval frame}.

The next result shows that if the integer translates of the generators of
a SIS form a frame, then the set of integer translations of the ``cutoffs''
of the generators remains a frame.

\begin{theorem}
Assume that $S$ is a SIS that is $\frac1n\Z$-invariant,
and that $\F \subseteq S$ is such that $\Tc_\Z(\F)$ is a frame for~$S$
with frame bounds $A$, $B$.
Then
$$\Tc_\Z(P_k\F)
\EQ \bigset{T_j f^k : f \in \F, \, j \in \Z}$$
is a frame for $U_k = \Sf(P_k\F)$ with frame bounds $A$, $B$.
Further,
$$\Tc_\Z\biggparen{\medcup_{k=0}^{n-1} P_k\F}
\EQ \bigset{T_j f^k : f \in \F, \, j \in \Z, \, k=0,\dots,n-1}$$
is a frame for~$S$ with frame bounds $A$, $B$.
\end{theorem}
\begin{proof}

By hypothesis,
\begin{equation} \label{transframe_def}
\forall\, g \in S, \quad
A \, \norm{g}_2^2
\LE \sum_{j \in \Z} \sum_{f \in \F} \, |\ip{g}{T_j f}|^2
\LE B \, \norm{g}_2^2.
\end{equation}

Suppose that $g \in U_k$.
Then since $P_k$ commutes with integer translations, we have
\begin{align*}
\sum_{j \in \Z} \sum_{f \in \F} \, |\ip{g}{T_j P_k f}|^2
& \EQ \sum_{j \in \Z} \sum_{f \in \F} \, |\ip{g}{P_k T_j f}|^2 \\[1 \jot]
& \EQ \sum_{j \in \Z} \sum_{f \in \F} \, |\ip{P_k g}{T_j f}|^2 \\[1 \jot]
& \EQ \sum_{j \in \Z} \sum_{f \in \F} \, |\ip{g}{T_j f}|^2.
\end{align*}
Combining this with~\eqref{transframe_def},
we see that $\Tc_\Z(P_k\F)$ is a frame for $U_k$ with frame bounds $A$,~$B$.

Suppose now that $g \in S$.
Then since $S$ is the orthogonal direct sum of the $U_k$, we have that
\begin{align*}
\sum_{j \in \Z} \sum_{f \in \F} \sum_{k=0}^{n-1} \, |\ip{g}{T_j P_k f}|^2
\EQ \sum_{k=0}^{n-1} \sum_{j \in \Z} \sum_{f \in \F} \, |\ip{P_k g}{T_j f}|^2
\LE B \sum_{k=0}^{n-1} \, \norm{P_k g}_2^2
\EQ B \, \norm{g}_2^2.
\end{align*}
The estimate from below is similar, so we see that
$\Tc_\Z\bigparen{\medcup_{k=0}^{n-1} P_k\F}$
is a frame for $S$, with frame bounds $A$,~$B$.
\end{proof}

\smallskip
\subsection{Characterization of $\frac1n$-invariance in terms of fibers}

A useful tool in the theory of shift-invariant spaces is based on
early work of Helson \cite{Hel64}.
An $L^2(\R)$ function is decomposed into ``fibers.''
This produces a characterization of SIS in terms of 
closed subspaces of $\ell^2(\Z)$ (the fiber spaces).
For a detailed description of this approach, see \cite{Bow00} and the
references therein.

\begin{definition}
Given $f \in L^2(\R)$ and $\omega \in [0,1)$,
the \emph{fiber} $\fhat_\omega$ of~$f$ at~$\omega$ is the sequence
$$\fhat_\omega \EQ \bigset{\fhat(\omega+k)}_{k \in \Z}.$$
\end{definition}

If $f$ is in $L^2(\R)$, then the fiber $\fhat_\omega$ belongs
to $\ell^2(\Z)$ for almost every $\omega \in [0,1)$.

\begin{definition}
Given a  subspace $V$ of $L^2(\R)$ and $\omega \in [0,1)$,
the \emph{fiber space} of~$V$ at~$\omega$ is
$$\Jc_V(\omega)
\EQ\overline{\bigset{\fhat_\omega :
    f \in V \text{ and } \fhat_\omega \in \ell^2(\Z)}},$$
where the closure is taken in the norm of $\ell^2(\Z)$.
\end{definition}
The map assigning to each $\omega$ the fiber space $\Jc_V(\omega)$
is known in the literature as the {\it range function} of $V$. 

For a proof that, for almost every~$\omega$, $\Jc_V(\omega)$ is a well-defined closed subspace
of $\ell_2(\Z)$  and that shift-invariant spaces can be characterized through range functions, see \cite{Bow00}, \cite{Hel64}.

We will need the following two results.

\begin{proposition}[\cite{Hel64}] \label{bownik-1}
If $S$ is a SIS, then
$$S \EQ \bigset{f \in L^2(\R) :
        \fhat_\omega \in \Jc_S(\omega) \text{ for a.e. } \omega}.$$
\end{proposition}

\begin{proposition} \label{ortho}
Let $S_1$ and $S_2$ be SISs.
If $S = S_1 \dotoplus S_2$, then
$$\Jc_S(\omega)
\EQ \Jc_{S_1}(\omega) \dotoplus \Jc_{S_2}(\omega), \quad\text{a.e. } \omega.$$
\end{proposition}

The converse of Proposition~\ref{ortho} is also true, but will not be needed.

Combining Theorem~\ref{main} with Proposition~\ref{bownik-1}
yields the following characterization of $\frac1n$-invariance
in terms of the fiber spaces. 

\begin{theorem}
Let $S$ be a SIS.
Then the following statements are equivalent.

\begin{enumerate}
\item[(a)]
$S$ is $\frac1n\Z$-invariant.

\smallskip
\item[(b)]
$\Jc_{U_k}(\omega) \subseteq \Jc_S(\omega)$
for almost every~$\omega$ and each $k=0,\dots,n-1$.
\end{enumerate}
\end{theorem}

For the finitely generated case we can obtain a slightly simpler
characterization of $\frac1n\Z$-invariance.

\begin{theorem} \label{fibers}
If $S$ is a finitely generated SIS, then
the following statements are equivalent.

\begin{enumerate}
\item[(a)]
$S$ is $\frac1n\Z$-invariant.

\smallskip
\item[(b)]
For almost every $\omega \in [0,1)$,
$$\Dim\bigparen{\Jc_S(\omega)}
\EQ \sum_{k=0}^{n-1} \Dim\bigparen{\Jc_{U_k}(\omega)}.$$
\end{enumerate}
\end{theorem}

\begin{proof} 
(a) $\Rightarrow$ (b).
If $S$ is $\frac1n\Z$-invariant then
$S = \dotbigoplus_{k=0}^{n-1} U_k$.
This is an orthogonal direct sum, so Proposition~\ref{ortho} implies that
$\Jc_S(\omega) = \dotbigoplus_{k=0}^{n-1} \Jc_{U_k}(\omega)$
for a.e.~$\omega$, with this sum also orthogonal.
The equality of dimensions in statement~(b) therefore holds.

\medskip
(b) $\Rightarrow$ (a).
Suppose that statement~(b) holds.
It is clear that the inclusion
$\Jc_S(\omega) \subseteq \dotbigoplus_{k=0}^{n-1} \Jc_{U_k}(\omega)$
holds for a.e.~$\omega$. 
Since the spaces $U_k$ are orthogonal,
Proposition~\ref{ortho} implies that,
for a.e.~$\omega$,
the spaces $\Jc_{U_k}(\omega)$ are also orthogonal.
Counting dimensions and applying statement~(b), we conclude that
$\Jc_S(\omega) = \dotbigoplus_{k=0}^{n-1} \Jc_{U_k}(\omega)$ for a.e.~$\omega$.

Suppose now that $f \in U_k$.
Then $\fhat_\omega \in \Jc_{U_k}(\omega) \subseteq \Jc_S(\omega)$
for a.e.~$\omega$, so Proposition~\ref{bownik-1} implies that $f \in S$.
Thus $U_k \subseteq S$ for each $k$, so it follows from
Theorem~\ref{main} that $S$ is $\frac1n\Z$-invariant.
\end{proof}
\begin{remark}
Given a SIS $V$, the function $D_V(\omega)\equiv \Dim\bigparen{\Jc_V(\omega)}$ defined 
for $\omega \in [0,1)$ is known in the literature  as the {\it Dimension function} or {\it Multiplicity function}
of  the shift-invariant space $V$. So, condition (b) of Theorem \ref{fibers} is a statements about dimension
functions of the shift invariant spaces involved.
\end{remark}

\medskip
\subsection{The Bownik decomposition and $\frac1n$-invariance}

In \cite{Bow00}, Bownik obtained a decomposition for general
shift-invariant spaces, extending the earlier works
\cite{dBDR94a} and \cite{dBDR94b},
which applied to the finitely generated case.
We will apply this decomposition to shift-invariant spaces
that are $\frac1n\Z$-invariant.

\begin{theorem}[Bownik] \label{Bownik}
Let $S \subseteq L^2(\R)$ be a SIS.
Then for each $j \in \N$ we can find a function
$\varphi_j \in L^2(\R)$ such that
$\Tc_\Z(\varphi_j)$ is a Parseval frame for $\Sf(\varphi_j)$,
and furthermore
$$S \EQ \dotbigoplus_{j \in \N} \Sf(\varphi_j).$$
\end{theorem}

Note that a consequence of this theorem is that every SIS always has 
a set of generators whose integer translates form a Parseval frame of the SIS.

By applying Theorem~\ref{Bownik} to each space~$U_k$,
we obtain the following result.

\begin{theorem}
Let $S$ be a SIS that is $\frac1n\Z$-invariant.
Then there exist functions $\varphi_{k,j} \in L^2(\R)$ such that
$$S \EQ \dotbigoplus_{k=0}^{n-1} \dotbigoplus_{j \in \N} \Sf(\varphi_{k,j}),$$
with the following properties holding.

\begin{enumerate}

\item[(a)] \label{p2}
$\Tc_\Z(\varphi_{k,j})$ is a Parseval frame for $\Sf(\varphi_{k,j})$.

\smallskip
\item[(b)] \label{p3}
$\Sf(\varphi_{k,j}) \subseteq U_k$ for each $j \in \N$, and
$$U_k \EQ \dotbigoplus_{j \in \N} \Sf(\varphi_{k,j}).$$

\item[(c)] \label{p4}
Each space $\Sf(\varphi_{k,j})$ is $\frac1n\Z$-invariant.
\end{enumerate}
\end{theorem}

\section{Finitely Generated Shift-Invariant Spaces and $\frac1n$-Invariance}

In this section we will apply some of the general results obtained so far
to the particular case of finitely generated shift-invariant spaces.
We will use the concept of the Gramian.
This is a  common tool in the study of finitely generated shift-invariant spaces; see for example \cite{dBDR94b},\cite{RS95},\cite{ACHM07}.

\smallskip
\subsection{Characterization of $\frac1n$-invariance in terms of the Gramian}

\begin{definition}
Let $\Phi= \set{\varphi_1,\dots,\varphi_m}$ be a collection of finitely
many functions in $L^2(\R)$.
Then the \emph{Gramian} $G_\Phi$ of $\Phi$ is the
$m \times m$ matrix of $\Z$-periodic functions 
\begin{equation} \label{gram}
[G_{\Phi}(\omega)]_{ij}
\EQ \Bigip{(\phihat_i)_\omega}{(\phihat_j)_\omega}
\EQ \sum_{k \in \Z} \phihat_i(\omega+k) \, \overline{\phihat_j(\omega+k)},
\qquad \omega\in \R,
\end{equation}
where $(\phihat_j)_\omega$ is the fiber of $\varphi_j$ at $\omega$.
\end{definition}

We consider now the SIS $S = \Sf(\Phi)$ generated by the set
$\Phi= \set{\varphi_1,\dots,\varphi_m}$.
It is known \cite{dBDR94b} that if $f \in \Sf(\Phi)$,
then there exist $\Z$-periodic functions $a_1,\dots,a_m$ such that
$$\fhat(\omega)
\EQ \sum_{j=1}^m a_j(\omega) \, \phihat_j(\omega),
\qquad\text{a.e. } \omega.$$
This implies that the fiber spaces $\Jc_S(\omega)$ are generated
by the fibers of the generators of $S$ at $\omega$
(see also \cite{Bow00}).
That is, for almost every~$\omega$ we have that
$$\Jc_S(\omega)
\EQ \Span\bigset{(\phihat_j)_\omega : j = 1,\dots,m}.$$
Therefore
$$\Dim\bigparen{\Jc_S(\omega)}
\EQ \rank[G_{\Phi}(\omega)]$$
for almost every~$\omega$.

In the same way, since the SIS $U_k$ is generated by
$\Phi^k = P_k\Phi = \set{\varphi_1^k, \dots, \varphi_m^k}$,
where $\varphi_j^k = P_k \varphi_j$,
we have for almost every~$\omega$ that the fiber spaces
$\Jc_{U_k}(\omega)$ satisfy
$$\Jc_{U_k}(\omega)
\EQ \Span\Bigset{\Bigparen{\widehat{\varphi_j^k}}_\omega :
    j = 1, \dots, m}.$$
Let us denote by $G_{\Phi^k}$ the Gramian matrix associated with
the generators of~$U_k$.
Then, as above we have that
$\Dim(\Jc_{U_k}(\omega)) = \rank[G_{\Phi^k}(\omega)]$
for almost every~$\omega$ and $k = 0, \dots, n-1$.
Now Theorem~\ref{fibers} can be restated in the following way.

\begin{theorem} \label{FSIS}
If $S = \Sf(\Phi)$ is the SIS generated by
$\Phi = \set{\varphi_1,\dots,\varphi_m}$,
then the following statements are equivalent.

\begin{enumerate}
\item[(a)]
$S$ is $\frac1n\Z$-invariant.

\smallskip
\item[(b)] For almost every~$\omega \in [0,1)$ we have
$$\rank[G_{\Phi}(\omega)] \EQ \sum_{k=0}^{n-1} \rank[G_{\Phi^k}(\omega)].$$
\end{enumerate}
\end{theorem}

\smallskip
\subsection{Implications for frequency support}

As a consequence of Theorem~\ref{FSIS} we deduce an interesting result
about the supports of the Fourier transforms of the generators of a SIS.

\begin{theorem} \label{supp-fhat}
Let $S = \Sf(\Phi)$ be the SIS generated by
$\Phi = \set{\varphi_1,\dots,\varphi_m}$, and define

$$E_j \EQ \bigset{\omega \in [0,1): \rank[G_{\Phi}(\omega)] = j},
\quad j = 0,\dots,m.$$
If $S$ is $\frac1n\Z$-invariant, then
for each interval $I \subseteq \R$ of length~$n,$ we have that
for each $h = 1, \dots, m$ 
$$\bigabs{\bigset{\omega \in I : \widehat{\varphi_h}(\omega) = 0}}
\GE \sum_{j=0}^{n-1} \, (n-j) \, |E_j|.$$
In particular if $n>m$ we have,
$$\bigabs{\bigset{\omega \in I : \widehat{\varphi_h}(\omega) = 0}}
\GE n-m.$$
\end{theorem}

\begin{proof}
The measurability of the sets $E_j$ follows from the results
of Helson \cite{Hel64}, e.g., see \cite{BK06} for an argument of this type.

It is enough to prove the theorem for the interval $I = [0,n)$.
 
We note that the set
$$K_1
\EQ \set{\omega \in [0,1) :
    (\widehat{\varphi_h})_\omega \in \ell^2(\Z) \text{ for } h = 1, \dots, m}$$
has full measure.
Therefore
$${K_n} \EQ \medcup_{k=0}^{n-1} (K_1 + k)$$
is a subset of $[0,n)$ of measure~$n$.

Fix any particular $j \in \set{0, \dots, m}$.
By Theorem~\ref{FSIS},
\begin{equation} \label{gram_eq}
\rank[G_{\Phi}(\omega)] \EQ \sum_{k=0}^{n-1} \rank[G_{\Phi^k}(\omega)],
\qquad\text{a.e.\ } \omega.
\end{equation}
Therefore, for $j  < n,$ if $\omega \in E_j$ then at least $n-j$ terms on the
right-hand side of equation~\eqref{gram_eq} must vanish.
That is, given such an~$\omega$, there exists a subset of
$\set{0, \dots, n-1}$
with at least $n-j$ elements such that for each $k$ in this subset
$$\rank[G_{\Phi^{k}}(\omega)] = 0.$$
In particular, for each $h \in \{1,\dots,m\}$ we have that, for $j  < n$ and $ \omega \in E_j,$ 
\begin{equation}\label{zero}
\# \bigset{k \in \{0, \dots, n-1\}: \widehat{\varphi_h}(\omega+k) = 0}\geq n-j.
 \end{equation}

On the other hand, using that the sets $E_j$ are disjoint, we have
\begin{eqnarray*}
\lefteqn{\bigset{\omega\in K_n : \widehat{\varphi_h}(\omega) = 0}}\\
 & \EQ &  \bigset{\omega+ k :\omega \in K_1,\; 0 \leq k \leq n-1 \text{ and }
\; \widehat{\varphi_h}(\omega+k) = 0}\\
& \EQ & \bigcup_{j=0}^m  \bigset{\omega+ k :\omega \in K_1 \cap E_j, \; 0 \leq k \leq n-1 \text{ and }
\; \widehat{\varphi_h}(\omega+k) = 0}.
\end{eqnarray*}

Consequently,  from (\ref{zero}) and the last equation it follows that,
\begin{eqnarray*}
\lefteqn{\bigabs{\bigset{\omega\in K_n : \widehat{\varphi_h}(\omega) = 0}}}\\
& \EQ &\sum_{j=0}^m \bigabs{ \bigset{\omega+ k :\omega \in K_1 \cap E_j,  0 \leq k \leq n-1 \text{ and }
 \widehat{\varphi_h}(\omega+k) = 0}}\\
  &  \EQ & \sum_{j=0}^m \int_{E_j} \# \bigset{k \in \{0, \dots, n-1\}: \widehat{\varphi_h}(\omega+k) = 0}\;d\omega \\
&  \geq &\sum_{j=0}^{n-1} (n-j) |E_j|.
\end{eqnarray*}

Furthermore, if $n > m$, 
\begin{equation*}
\sum_{j=0}^m \, (n-j) \, |E_j|
\GE \sum_{j=0}^m \, (n-m) \, |E_j|
\EQ (n-m).
\end{equation*}
The measurability of the function
 $\omega \longmapsto \# \bigset{k \in \{0, \dots, n-1\}: \widehat{\varphi_h}(\omega+k) = 0},$
 follows from the fact that
 \begin{equation*}
\begin{split} &\bigset{w \in E_j:   \# \bigset{k \in \{0, \dots, n-1 \}: \widehat{\varphi_h}(\omega+k) = 0} \geq s} \\
 \EQ  &\bigcup_{0\leq k_1<\dots<k_s < n}\bigcap_{i=1}^s \bigset{\omega \in E_j:
 \widehat{\varphi_h}(\omega+k_i) = 0}.
 \end{split}
 \end{equation*}

\end{proof}
 
Note that if $S$ is a principal SIS, say $S = \Sf(\varphi)$,
then the Gramian is scalar-valued; since
$$G_\varphi(\omega)
\EQ \bigip{\phihat_\omega}{\phihat_\omega}
\EQ \sum_{k \in \Z} \, |\phihat(\omega+k)|^2.$$
Applying Theorem~\ref{supp-fhat} to this case,
we obtain the following corollary.

\begin{corollary} \label{measure-zero}
Let $\varphi \in L^2(\R)$ be given.
If the SIS $\Sf(\varphi)$ is $\frac1n\Z$-invariant for some $n > 1$,
then $\phihat$ must vanish on a set of infinite Lebesgue measure.
Furthermore, for each interval $I \subseteq \R$ of length~$n$, we have that 
$$\bigabs{\bigset{\omega \in I : \phihat(\omega) = 0}}
\GE n \, |E_0| + (n-1) \, |E_1|
\GE n-1,$$
where
$E_0 = \set{\omega \in [0,1): G_{\varphi}(\omega) = 0}$ and $E_1 = \set{\omega \in [0,1): G_{\varphi}(\omega) \neq 0}$
\end{corollary} 
 
This yields the following fact regarding the order of invariance of a
principal SIS generated by a compactly supported function.
 
\begin {proposition} \label{compact}
If a nonzero function $\varphi \in L^2(\R)$ has compact support,
then $\Sf(\varphi)$ has invariance order one.
That is $\Sf(\varphi)$ is not $\frac1n\Z$-invariant for any $n > 1$.
\end{proposition}
 
\begin{proof}
Because of Corollary~\ref{measure-zero},
if $\Sf(\varphi)$ is $\frac1n\Z$-invariant with $n>1$,
then $\phihat$ must vanish on a set of positive measure.
Since $\varphi$ has compact support, the Paley--Wiener Theorem
implies that $\varphi = 0$~a.e.
\end{proof}

It is not difficult to construct a function $\varphi \in L^2(\R)$ such that
$\phihat$ is compactly supported in frequency yet the SIS $\Sf(\varphi)$
is not translation-invariant.
In fact, we have the following consequence of Corollary~\ref{measure-zero}.

\begin{corollary} \label{TI}
If $\varphi \in L^2(\R)$ and $\Sf(\varphi)$ is translation-invariant,
then $|\supp(\phihat)| \le 1$.
\end{corollary} 

\begin{remark}\label{remark}
As one of the referees pointed out, Proposition \ref{compact} is known and follows readily from Proposition \ref{BDR2}.

Likewise, Corollary \ref{TI} can  be obtained from properties of the dimension function of a SIS.
For this, observe  that since $\Sf(\varphi)$ is translation invariant then,

$$D_{\Sf(\varphi)} = \sum_{k\in\Z}\chi_{supp(\hat\varphi)}(\omega+k).$$

Now using that $\Sf(\varphi)$ is principal, we have $D_{\Sf(\varphi)}\leq 1.$
Thus integrating both sides over $[0,1],$ yields $|\supp(\phihat)| \le 1$.
For properties of the dimension function of a SIS, see for example \cite{BM99} or \cite{BR03}.
\end{remark}

\begin{remark}
Assume now that $\varphi \in L^2(\R)$ and $\Sf(\varphi)$ is translation-invariant.
Using the argument in Remark \ref{remark} we have that
 $ \sum_{k\in\Z}\chi_{\supp(\hat\varphi)}(\omega+k) \leq 1.$
This implies that supp$(\hat\varphi)$ is a subset of a set of representatives of the quotient $\R/\Z$. 
That is, supp$(\hat\varphi)$ is a subset of a tile of $\R,$ what is a refinement of Corollary \ref{TI}.
Translation-invariance of multiresolution analyses is connected to tilings of $\R$. See for example \cite{Mad92}.
\end{remark}
\smallskip
\subsection{Application to multiresolution analyses}

For definitions and details on wavelets
and multiresolution analyses, see \cite{Dau92} or \cite{Mal89}.

Suppose that $\set{V_j}_{j \in \Z}$ is a multiresolution analysis (MRA)
of $L^2(\R)$.
By definition, $V_0 = \Sf(\varphi)$ for some $\varphi \in L^2(\R)$,
called the \emph{scaling function},
and $V_j$ is the image of $V_0$ under the unitary operator
$D_{2^j}f(x) = 2^{j/2} \, f(2^j x)$.

The preceding results imply that if $\varphi$ is compactly supported
(as is the case for the Daubechies scaling functions, for example),
then the SIS $V_0$ has invariance order exactly~$1$.
Thus $V_0$ is invariant only under integer translations.
The same remarks apply to the associated wavelet $\psi$ and
wavelet SIS $W_0 = \Sf(\psi)$ if $\psi$ is compactly supported.

Further, if $\varphi$ is compactly supported, then at resolution level~$j$,
the subspace $V_j$ is invariant exactly under translations $2^j\Z$,
and similarly for the wavelet space $W_j$ if $\psi$ is compactly supported.
This includes all the spaces associated with the Daubechies scaling
functions and wavelets, for example.

\section{Higher dimensions}

The results of this article are for the line. The higher-dimensional case is much more
involved due to the more complex structure of the closed additive subgroups of $\R^d.$ 
This will be the subject of a forthcoming article.

\section{Acknowledgements}
We thank Magal\'\i\ Anastasio for carefully reading the manuscript.
We also thank the referees for very interesting observations that
contributed to improvement of the paper and the argument in Remark \ref{remark}.

\end{document}